\documentclass{article}
\usepackage{amssymb,amsmath,amsfonts,epsfig,latexsym}

\def\Z{\ensuremath{\mathbb{Z}}} 
\def\P{\ensuremath{\mathbb{P}}}
\def\R{\ensuremath{\mathbb{R}}}

\def\OO{\ensuremath{\mathcal{O}}}

\def\<{\ensuremath{\langle}}
\def\>{\ensuremath{\rangle}}

\begin{document}
    
\title{A smooth, complete threefold with no nontrivial nef line bundles}           
    
\author{Sam Payne}

\maketitle

\begin{abstract}
   We present a smooth, complete toric threefold with no nontrivial
   nef line bundles.  This is a counterexample to a recent conjecture
   of Fujino.
\end{abstract}

Let $X$ be a complete algebraic variety.  A line bundle $L$ on $X$ is
said to be nef if $L \cdot C$, the degree of $L$ restricted to $C$, is
nonnegative for every curve $C \subset X$.  If $X$ has no nontrivial
nef line bundles, then $X$ admits no nonconstant maps to projective
varieties.  Indeed, if $f:X \rightarrow \P^{n}$ is a nonconstant
morphism, then $f^{*} \OO(1)$ is a nef line bundle that is not
isomorphic to $\OO_{X}$.  One example of a complete variety with no
nontrivial nef line bundles, due to Fulton, is the singular toric
threefold corresponding to the fan over a cube with one ray displaced,
which has no nontrivial line bundles at all \cite[pp.25-25,
72]{Fulton}.  The literature also contains a number of examples of
Moishezon spaces with no nontrivial nef line bundles.  See
\cite[3.3]{Nakamura}, \cite[5.3.14]{Kollar}, \cite{Oguiso}, and
\cite{Bonavero}.  Recently, Fujino has given examples of complete,
singular toric varieties with arbitrary Picard number that have no
nontrivial nef line bundles.  Fujino conjectured that $NE(X)$ is
properly contained in $N_{1}(X)$ for any smooth, complete toric
variety $X$.  This conjecture is equivalent to the statement that
every smooth, complete toric variety has a nontrivial nef line bundle
\cite[Remark 3.1, Conjecture 4.5]{Fujino}.  The purpose of this note
is to give a counterexample to Fujino's conjecture; we present a
smooth, complete toric threefold with no nontrivial nef line bundles.

\vspace{10 pt}

Let $\Sigma$ be the fan in $\R^{3}$ whose rays are generated by
\[
\begin{array}{llll}
    v_{1} = (1,1,1), & v_{2} = (-1,1,1), & v_{3} = (-1,-1,1), & v_{4} =
    (1,-1,1), \\ v_{5} = (1,0,-1), & v_{6} = (0,1,-1), & v_{7} = (-1, 0,
    -1), & v_{8} = (0,-1,-1), 
\end{array}
\]
and whose maximal cones are $\<v_{1}, v_{2}, v_{3}, v_{4} \>$, $\<
v_{5}, v_{6}, v_{7}, v_{8} \>$ and
\[
\begin{array}{llll}
    \< v_{1}, v_{2}, v_{5} \>, & \< v_{2}, v_{3}, v_{6} \>, & 
    \< v_{3}, v_{4}, v_{7} \>, & \< v_{4}, v_{1}, v_{8} \>, \\
    \< v_{2}, v_{5}, v_{6} \>, & \< v_{3}, v_{6}, v_{7} \>, &
    \< v_{4}, v_{7}, v_{8} \>, & \< v_{1}, v_{8}, v_{5} \>.
\end{array}
\]
The nonzero cones of $\Sigma$ are the cones over the faces of the 
following nonconvex polyhedron, which has a rotational symmetry about 
the vertical axis of order 4.

\hspace{50 pt}{\includegraphics{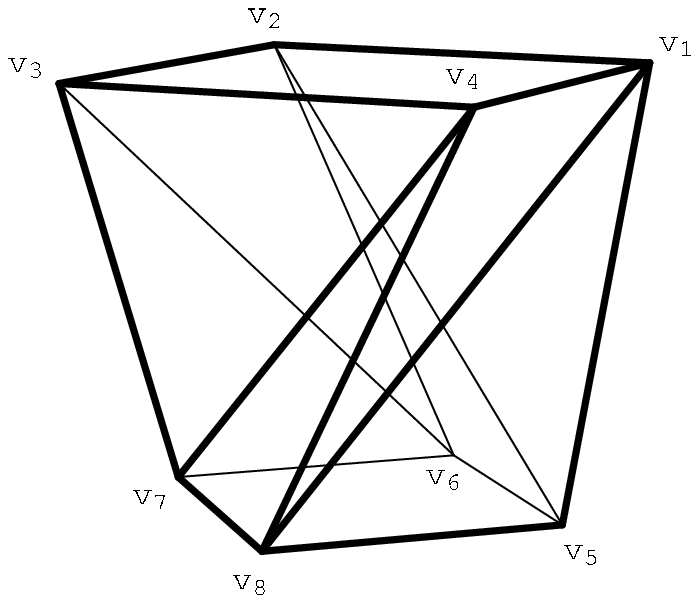}}

Let $\Delta$ be the fan obtained from $\Sigma$ by successive stellar
subdivisions along the rays spanned by $v_{9} = (0,0,-1),$ $v_{10} =
(0,0,1)$, and
\[
\begin{array}{llll}
    v_{11} = (1,0,1), & v_{12} = (0,1,1), 
    & v_{13} = (-1,0,1), & v_{14} = (0,-1,1).
\end{array}
\]
We claim that $X = X(\Delta)$, the toric threefold corresponding to the
fan $\Delta$ with respect to the lattice $\Z^{3} \subset \R^{3}$, is 
smooth and complete and that $X$ has no nontrivial nef line bundles. 
It is easy to check that $|\Delta| = \R^{3}$, so $X$ is complete, and
that each of the 24 maximal cones of $\Delta$ is spanned by a basis
for $\Z^{3}$, so $X$ is smooth.  It remains to show that $X$ has no
nontrivial nef line bundles.

Suppose $D = \sum d_{i} D_{i}$ is a nef Cartier divisor, where $D_{i}$
is the prime $T$-invariant divisor corresponding to the ray spanned by
$v_{i}$.  Since nef line bundles on toric varieties are globally
generated, $\OO(D)$ has nonzero global sections.  Therefore we may
assume that $D$ is effective, i.e.\ each $d_{i} \geq 0$.  We will show
that $D = 0$.  Note that $\< v_{2}, v_{5}, v_{6} \>$ is a cone in
$\Delta$, and
\[
   v_{1} \, = \, 3v_{2} + 4 v_{5} -2 v_{6}.
\]
It then follows from the convexity of the piecewise linear function 
$\Psi_{D}$ associated to $D$ that
\[
   d_{1} + 2d_{6} \ \ \geq \ \ 3d_{2} + 4 d_{5}.
\]   
Similarly,
\[
\begin{array}{lll}
   d_{2} + 2d_{7} & \geq & 3d_{3} + 4 d_{6}, \\  
   d_{3} + 2d_{8} & \geq & 3d_{4} + 4 d_{7}, \\
   d_{4} + 2d_{5} & \geq & 3d_{1} + 4 d_{8}.
\end{array}   
\]
Adding the four above inequalities, we have
\begin{eqnarray*}
    \lefteqn{ d_{1} + d_{2} + d_{3} + d_{4} + 2(d_{5} + d_{6} + 
    d_{7} + d_{8}) \ \geq } \\
    & & 3(d_{1} + d_{2} + d_{3} + d_{4}) + 4(d_{5} + d_{6} + d_{7} + d_{8}).
\end{eqnarray*}    
Since all of the $d_{i}$ are nonnegative, it follows that $d_{i} = 0$
for $1 \leq i \leq 8$.  Since $v_{1}, \ldots, v_{8}$ positively span
$\R^{3}$, we have $P_{D} = \{ 0 \}$, and hence $D = 0$.

\vspace{10 pt}

The above argument shows that any complete fan in $\R^{3}$ containing
the cones $\< v_{2}, v_{5}, v_{6} \>,$ $\<v_{3}, v_{6}, v_{7} \>$,
$\<v_{4}, v_{7}, v_{8} \>$, and $\< v_{1}, v_{5}, v_{8} \>$
corresponds to a complete toric variety with no nef line bundles.  In
particular, stellar subdivisions of the remaining cones in $\Delta$
lead to examples of smooth, complete toric threefolds with arbitrarily
large Picard number and no nontrivial nef line bundles.

\end{document}